\documentclass{amsart}

\usepackage[centertags]{amsmath}
\usepackage{amsfonts}
\usepackage{amssymb}
\usepackage{amsthm}
\usepackage{graphicx}

\begin{document}

\title{Polynomial ideals and directed graphs}
\author{GIUSEPPA CARRA' FERRO and DANIELA FERRARELLO}
\address{Department of Mathematics and Computer Science\\University of Catania\\
Viale Andrea Doria 6, 95125 Catania, Italy}
\email {carra@dmi.unict.it\\
ferrarello@dmi.unict.it}
\date{}

\newtheorem{defn}{DEFINITION}%[section]
\newtheorem{thm}{THEOREM}%[section]
\newtheorem{ex}{EXAMPLE}%[section]
\newtheorem{rem}{REMARK}%[section]
\newtheorem{prop}{PROPOSITION}%[section]
\newtheorem{cor}{COROLLARY}%[section]

\renewcommand{\qed}{\hfill$\Box$\smallskip}
\newpage

%\medskip
%\noindent
% Short Title: Polynomial ideals and directed graphs

% \medskip
% \noindent
% Address of the author, that has to correct the proofs:\\
% Daniela Ferrarello, ferrarello@dmi.unict.it

% \newpage
\begin{abstract}
In this paper it is shown that it is possible to associate several
polynomial ideals to a directed graph $D$  in order to find
properties of it. In fact by using algebraic tools it is possible
to give appropriate procedures for automatic reasoning on cycles
and directed cycles of graphs.
\end{abstract}
\maketitle

{\bf Keywords}: ideals, digraphs, algorithms.

\section{Introduction}
The aim of this paper is to study some  problems in graph theory
with commutative algebra tools.  Connections between simplicial
complexes and polynomial
 rings were first studied in ( \cite{Sta96}).  There were also studies on the
 connections between ideals and undirected graphs, as in Simis, Vasconcelos and Villarreal
 (\cite{SVV94} and \cite{SVV98}). In 1995 Villarreal (\cite{Vi95}) and in 1998
 Hibi and Ohsugi (\cite{HO98}) associated a
 toric ideal, called the {\em edge ideal}, to an undirected graph. They found some relations
  between ideal properties and even closed walks of the associated graph.
Other studies of undirected graphs with the edge ideal can be found
in  Villarreal  (\cite{Vi01}) and  Herzog and Hibi (\cite{HH05}). By
using these papers the authors found other ideals associated to a
graph and other properties on cycles and
minimal vertex covers (\cite{CF04}).\\
In this paper we extend the paper \cite{CF04} to the case of
directed graphs by founding a toric ideal, such that its
generators are in one to one correspondence with directed and
undirected cycles. The existence of such ideal was proved in
different way in 2005 by Reyes (\cite{Re05}) and in 2006 by
Gitler, Reyes and Villarreal (\cite{GRV06}. Some binomial ideals
associated to a digraph can be also found in \cite{IsIm00} and
\cite{BPS01}, but there is no characterization of cycles.
Furthermore we introduce the notion of sink and source covers and
we find characteristic conditions for directly bipartite graphs.
Many properties of a digraph $D$ are studied trough some
corresponding properties of two undirected graphs associated to
$D$. Relations between Cohen-Macaulay property of such undirected
graphs and properties of the digraph will be investigated in
another paper.\\
Digraphs are very useful for many applications
like in computational molecular biology (\cite{dJ02}), in the
study of phylogenetic trees (\cite{Er04} and \cite{SS05}) and in
the minimum cost flow problem in networks, which has many physical
applications (\cite{GT89}). By using the packages {\em Groebner}
and {\em networks} of Maple 10 we have procedures for automatic
deduction in graph theory in order to find bases of directed and
undirected cycles of the given digraph and to check the existence
of sink and source covers.

\section{Preliminary tools} \label{tools}
\subsection{Gr\"{o}bner Bases}
In this section we introduce some basic notions and properties of
toric ideals. Let ${\bf N}_{0}$=$\{ 0,1,2,\ldots,n,\ldots \}$ and
let $X_{1},\ldots,X_{n}$ be $n$ variables. Let $K$ be a field of
characteristic zero. Let $A=K[X_{1},\ldots,X_{n}]$ and let
$PP(X_{1},\ldots,X_{n})$ =\{$X_{1}^{a_{1}} \cdots X_{n}^{a_{n}}$:
$(a_{1},\ldots,a_{n}) \in {\bf N}_{0}^{n}$\} be the set of power
products in $\{X_1,\ldots,X_n\}$, that is equal to the set $T_{A}$
of terms of $A$.

\begin{defn}
A term ordering $\sigma$ on $T_{A}$  is a total order such that:\\
(i) $1<_{\sigma} t$ for all $t \in T_{A} \setminus \{1\}$; \\
(ii) $t_{1}<_{\sigma} t_{2}$ implies $t_{1}t'<_{\sigma} t_{2}t'$
for all $t' \in T_{A}$.
\end{defn}

\smallskip
\noindent If $\sigma$ is a term ordering on $T_{A}$ and $f \in A$,
then $M_{\sigma}(f)$ is the monomial $c_{j}t_{j}$ iff
$t_{i}<_{\sigma}t_{j}$ for all $i \neq j$, $i=1,\ldots,r$.
$M_{\sigma}(f)$ is the {\em leading monomial} of f.
$t_j$=$T_{\sigma}(f)$ is the {\em leading term} of $f$, while
$M_{\sigma}(f)$ is the {\em leading monomial} of f.

\begin{defn}
Let $f$ be a polynomial in $A$, let $F=(f_{1},\ldots,f_{r})$ be a
finite subset of $A$ and let $\sigma$ be a term ordering on $T_A$.
If $f$=$\sum_{i=1,\ldots,s}c_{i}t_{i} \in A$, then  $f$ is {\em
reduced with respect to} $F$ iff $t_{i} \neq tT_{\sigma}(f_h)$ for
all $t \in T_{A}$, all $i=1,\ldots,s$ and all $h=1,\ldots,r$.
\end{defn}

\begin{defn}
Let $\sigma$ be a term ordering on $T_A$ and let $I$ be an ideal
in $A$. The monomial ideal $M_{\sigma}(I)=(M_{\sigma}(f): f \in
I)$ is the {\em initial ideal} of $I$.
\end{defn}

\begin{defn}
{\em (\cite{Bu76})} Let $I$ be an ideal in $A$ and let $\sigma$ be
a term ordering on $T_A$. If $I=(f_{1},\ldots,f_{r})$, then
$\{f_{1},\ldots,f_{r}\}$ is a Gr\"{o}bner basis of $I$ with
respect to $\sigma$ on $T_{A}$ iff
$M_{\sigma}(I)=(M_{\sigma}(f_{1}),\ldots,M_{\sigma}(f_{r}))$.\\
$\{f_{1},\ldots,f_{r}\}$ is a reduced Gr\"{o}bner basis iff
$M_{\sigma}(f_{h})$=$T_{\sigma}(f_{h})$  and $f_h$ is reduced with
respect to $F \setminus f_{h}$  for all $h=1,\ldots,r$.
\end{defn}

\smallskip
\noindent Every ideal $I$ in $A$ has a Gr\"{o}bner basis and a
reduced Gr\"{o}bner basis with respect to a given term ordering
$\sigma$ (\cite{Bu76}).

\begin{defn}\label{univ}
Let $I$ be a nonzero ideal in $A$. The {\em universal Gr\"{o}bner
basis of} $I$ is the union of all reduced Gr\"{o}bner bases of
$I$.
\end{defn}

\subsection{Toric Ideals}
Here we introduce the notion and  some properties of a toric
ideal. Let $K$ be a field and let $A=K[X_{1},\ldots,X_{n}]$ be as
above. let $B=K[X_{1},\ldots,X_{n},X_{1}^{-1},\ldots,X_{n}^{-1}]$
be the Laurent polynomial ring in the indeterminates
$\{X_{1},\ldots,X_{n}\}$ and let $EPP(X_{1},\ldots,X_{n})$
=\{$X_{1}^{a_{1}} \cdots X_{n}^{a_{n}}$: $(a_{1},\ldots,a_{n}) \in
{\bf Z}^{n}$\} be the set of power products in
$\{X_1,\ldots,X_n,X_{1}^{-1},\ldots,X_{n}^{-1}\}$, that is equal
to the set $ET_{A}$ of extended terms of $A$.

\begin{defn}\label{deftoric}
{\em (\cite{St95})}. Let
$M$=$(m_{ij})_{i=1,\ldots,m,j=1,\ldots,n}$ be a $(m,n)$-matrix
with $m_{ij}$ in ${\bf Z}$ for all $i,j$. Let $\pi_{M}$ : ${\bf
N}_{0}^{n} \longrightarrow {\bf Z}^{m}$ be the semigroup
homomorphism defined by
$\pi_{M}(u_1,\ldots,u_n)=(\sum_{j=1,\ldots,n}u_{j}m_{1j},\ldots,
\sum_{j=1,\ldots,n}u_{j}m_{mj})$. \\Let  $\exp_{n}$ : ${\bf
N}_{0}^{n} \longrightarrow PP(X_{1}, \ldots, X_{n})$ be the
semigroup isomorphism defined by  $\exp_{n}(u_1,\ldots,u_n)=
\prod_{j=1,\ldots,n}X_{j}^{u_{j}}$ and let  $\exp_{\bf Z}$ : ${\bf
Z}^{m} \longrightarrow EPP(t_{1},\ldots, t_{m})$ be the semigroup
isomorphism defined by $\exp_{\bf Z}(a_1,\ldots, a_m)$=
$\prod_{i=1,\ldots,m}t_{j}^{a_{i}}$. \\Let $\pi$ :
$PP(X_{1},\ldots,X_{n}) \longrightarrow EPP(t_{1},\ldots,t_{m})$
be the semigroup homomorphism, that is induced by $\pi_{M}$ and it
is defined by $\pi(X_j)=\prod_{i=1,\ldots,m}t_{i}^{m_{ij}}$ for
all $j=1,\ldots,n$. $\pi$ extends uniquely to the homomorphism of
semigroup algebras  $\pi':K[X_1,\ldots,X_n] \longrightarrow
K[t_1,\ldots,t_m,t_{1}^{-1},\ldots, t_{m}^{-1}]$ defined by
$\pi'(X_j)$ =$\pi(X_j)$ for all $j=1,\ldots,n$.
$I_{M}$=$ker(\pi')$ is an ideal in $A=K[X_1,\ldots,X_n]$, that is
called the {\em toric ideal of} $M$.
\end{defn}

\begin{rem}
Let $M$ be a matrix as above and let $I_M$ be the toric ideal of
$M$. It is not too hard to show that $I_M$ is an elimination
ideal. More precisely
$I_M$=$($$X_j-\prod_{i=1,\ldots,m}t_{i}^{m_{ij}},
t_{i}t_{i}^{-1}-1$: $j=1,\ldots,n$, $i=1,\ldots,m$$)$ $\cap$
$K[X_1,\ldots,X_n]$ $($\cite{St95}$)$.\\
 If we put $z_i=t_{i}^{-1}$ for all $i=1,\ldots,n$, then
$I_M$ is equal to the toric ideal $I_{M'}$, where
$M'$=$(m_{ij}')_{i=1,\ldots,2m,j=1,\ldots,n}$ is a $(2m,n)$-matrix
with $m_{ij}'$ in ${\bf N}$ for all $i,j$. Moreover
$m_{ij}'$=$m_{ij}$ for all $i=1,\ldots,m$ with $m_{ij} \in {\bf
Z}^{+}$, $m_{ij}'$=$0$ for all $i=1,\ldots,m$ with $m_{ij} \in
{\bf Z}^{-}$, $m_{ij}'$=$-m_{(i-m)j}$  for all $i=m+1,\ldots,2m$
with $m_{ij} \in {\bf Z}^{-}$ and $m_{ij}'$=$0$ for all
$i=m+1,\ldots,2m$ with $m_{ij} \in {\bf Z}^{+}$.
\end{rem}

\smallskip
\noindent It is well known that every toric ideal is a prime
binomial ideal. Other properties of binomial ideals can be found
in (\cite{ES96}), while properties of toric ideals can be found in
(\cite{St95}) and (\cite{KR05}, chap.2).

\begin{thm}
{\em (\cite{St95})}. $I_M$ is generated by binomials of the type
$X^{u^{+}}-X^{u^{-}}$, where $u^{+}, u^{-} \in {\bf Z}^{n}$ are
non negative with disjoint support.
\end{thm}

\begin{defn}
A binomial $X^{u^{+}}-X^{u^{-}}$ is {\em primitive} if there
exists no other binomial $X^{v^{+}}-X^{v^{-}} \in I_M$ such that
$X^{v^{+}}$ divides $X^{u^{+}}$ and $X^{v^{-}}$ divides
$X^{u^{-}}$. A binomial $X^{u^{+}}-X^{u^{-}}$ in $I_M$ is a {\em
circuit} if $supp(u^{+}-u^{-})$ is minimal with respect to
inclusion and the
coordinates of $u^{+}-u^{-}$ are relatively prime. \\
$Gr_{M}$=\{ primitive binomials in $I_M$\} is  the {\em Graver
basis of} $I_M$.
\end{defn}

\begin{thm}
{\em (\cite{St95})}. Let $U_M$ be the universal Gr\"{o}bner basis
of $I_M$ and let $C_M$ be the set of all circuits in
$I_M$. Then:\\
(i) every binomial in $U_M$  is primitive.\\
(ii) $C_M \subseteq U_M \subseteq Gr_{M}$ for every $M$.
\end{thm}

\subsection{Graphs and digraphs}

In this paper $G$=$(V(G),E(G))$ will be a finite graph with
$V(G)=\{v_1,\ldots,v_n\}$ and $E(G)=\{e_1,\ldots,e_m\}$.
Furthermore  $[v_{i},v_{j}]$ denotes the oriented edge from
$v_{i}$ to $v_{j}$, while
every non oriented edge between $v_{i}$ and $v_{j}$ is denoted by $\{v_i,v_j\}$.\\
The {\em underlying graph} $G_u$ of a directed graph $G$ is the
undirected graph with $V(G_u)=V(G)$ and the same undirected edges
of $G$. Often $D$ will denote a directed graph, that will be also
 called a {\em digraph}.
 All graphs in this paper will be
simple, i.e. without multiple edges.

%\smallskip
%\noindent
%Here we introduce some definitions from graph theory,
%that we shall use in the paper.

\begin{defn}
Let $D$ be a digraph. A vertex $v_i$ in $V(D)$ is called a {\em
source} if no edge is directed into $v_i$. A vertex $v_i$ in $D$
is called a {\em sink} if every adjacent edge is directed into
$v_i$.
\end{defn}

\begin{defn}
Let $D$ be a digraph. A {\em walk} of length $n$ from a vertex
$v_i$ to a vertex $v_j$ in $D$ is a sequence of vertices
$v_i$=$v_{i(1)}$, \ldots, $v_j$=$v_{i(n+1)}$, such that  either
$[v_{i(h)}v_{i(h+1)}]$ or $[v_{i(h+1)}v_{i(h)}]$ is in $E(D)$ for
all $h=1,\ldots,n$. A walk is called a {\em directed walk} if
$[v_{i(h)}v_{i(h+1)}] \in E(D)$ for all $h$. If
$v_{i(1)}=v_{i(n+1)}$ in a directed walk, then it is called a {\em
directed cycle}. A walk  is called {\em simple} if there are no
repeated edges, while it is called {\em elementary} whenever there
are no repeated vertices.
\end{defn}

\begin{defn}
Let $G$ be an undirected graph. $G$ is {\em bipartite} if its
vertices can be divided in two sets, such that no edge connects
vertices in the same set. Equivalently $G$ is {\em bipartite} iff
all cycles in
$G$ are even.\\
$G$ is {\em acyclic} if it has no cycle, while it is a {\em tree}
if it is connected and acyclic.
\end{defn}

\smallskip
\noindent
Some generalizations of such definitions in the directed
case are the following ones.

\begin{defn}
Let $D$ be a digraph. $D$ is called {\em directly bipartite}
whenever its underlying graph $D_u$ is bipartite with bipartition
sets $W$ and $W'$, such that every $w \in W$ is a source in $D$,
and every $w' \in W'$ is a sink in $D$. $D$ is called a {\em
directed acyclic graph}, {\em DAG} for short, when there are no
directed cycles.
\end{defn}

\section{Binomial ideals arising from a digraph}
Here we will introduce the binomial  ideals, associated to a
digraph, that we will use in the paper. First of all we introduce
some binomial ideals associated to the edges of a digraph.

\begin{defn}
The binomial {\em extended diedge ideal} of a digraph $D$ is the
ideal $I(D,E)$=$($ $e_h-z_{i}v_{j}, z_{i}v_{i}-1$:
$e_h=[v_{i},v_{j}] \in E(D)$, $i=1,\ldots,n$ $)$ in
$K[e_1,\ldots,e_m,v_1,\ldots,v_n,z_{1},\ldots,z_{n}]$. The ideal
$I(E)_D=I(D,E)\cap K[e_1,\ldots, e_m]$ is the {\em binomial diedge
ideal of the digraph} $D$.
\end{defn}

The definitions as above extend the analogous ones given in the case
of undirected graphs as in \cite{HO98} and \cite{CF04}

\begin{rem}
$I(E)_D$ is the toric ideal of the matrix $IM(D)^{t}$, that is the
transpose of the incidence matrix
$IM(D)=(a_{ih})_{i=1,\ldots,n,h=1,\ldots,m}$ of $D$ defined by
$a_{ih}=-1$ if $e_{h}$ leaves $v_{i}$, $a_{ih}=1$ if $e_{h}$
arrives to  $v_{i}$ and $a_{ih}=0$ if $v_{i} \notin e_{h}$ for
every $v_{i} \in V(E)$ and $e_{h} \in E(D)$.
\end{rem}

\begin{ex} \label{D1}
Let $D_1$ be the digraph with $V(D_1)=\{v_1, v_2, v_3, v_4,v_5\}$
and $E(D_1)$ = $\{e_1=[v_1,v_2]$, $e_2=[v_2,v_3]$,
$e_3=[v_3,v_1]$, $e_4=[v_1,v_4]$, $e_5=[v_3,v_4]$,
$e_6=[v_3,v_5]\}$. The binomial extended diedge ideal is
$I(D_1,E)=(e_1-z_1 v_2, e_2- z_2 v_3, e_3- z_3 v_1, e_4-z_1 v_4,
e_5-z_3 v_4, e_6-z_3 v_5, z_1 v_1-1, z_2v_2-1, z_3v_3-1, z_4v_4-1,
z_5v_5-1)$. The binomial diedge ideal of $D_1$ is
$I(E)_{D_1}$=$(e_1e_2e_3-1, e_3e_4-e_5, e_2 e_1 e_5-e_4)$.
\end{ex}

\begin{ex} \label{D2}
Let $D_2$ be the digraph with $V(D_2)=\{v_1,v_2,v_3,v_4,v_5\}$ and
$E(D_2)$ = $\{e_1=[v_1,v_2]$, $e_2=[v_2,v_3]$, $e_3=[v_1,v_3]$,
$e_4=[v_4,v_3]$, $e_5=[v_3,v_5]\}$. The binomial extended diedge
ideal of $D_2$ is $I(D_2,E)=(e_1-z_1 v_2, e_2- z_2 v_3, e_3- z_1
v_3, e_4-z_4 v_3, e_5-z_3 v_5, z_1 v_1-1, z_2v_2-1, z_3v_3-1,
z_4v_4-1, z_5v_5-1)$ while the binomial diedge ideal of $D_2$ is
$I(E)_{D_2}$=$(e_1e_2-e_3)$.
\end{ex}

\section{The undirected graph $H_{D}$ associated to $D$}
Here we associate an undirected bipartite graph $H_{D}$ to a
digraph $D$. The introduction of such graph allows to prove
properties of the cycles of $D$ trough the properties of the
cycles of $H_{D}$.

\begin{defn}
Let $D$ be a digraph.  Let $H_D$ be the undirected graph with
$V(H_D)$ = $V(D) \cup \{z_{1},\ldots,z_{n}\}$ and $E(H_D)$ =
\{$e=\{z_{i},v_{j}\}$: $\{v_{i},v_{j}\} \in E(D)$\} $\cup$
\{$f_{i}=\{z_{i},v_{i}\}$: $i=1,\ldots,n\}$. Let
$R=K[v_{1},z_{1},f_{1}, \ldots, v_{n},z_{n},f_{n}, e_{1}, \ldots,
e_{m}]$ and let $\pi: R \rightarrow R/(f_{1}-1,\ldots f_{n}-1)$ be
the canonical ring homomorphism defined by $\pi(v_{i})=v_{i}$,
$\pi(z_{i})=z_{i}$, $\pi(f_{i})=1$, $\pi(e_{j})=e_{j}$ for all
$i=1,\ldots,n$ and $j=1,\ldots,m$.
\end{defn}

\smallskip
\noindent If $I(H_D,E)$ is the binomial extended edge
ideal of the undirected graph $H_D$ as in \cite{CF04}, then
$\pi(I(H_D,E))=I(D,E)$ by definition of $\pi$.

\smallskip
\noindent
The graph $H_D$ has some properties, that are
independent on the properties of the graph $D$.

\begin{defn} A {\em matching} of an undirected graph $G$ is a
subset of independent edges of $G$(i.e. edges that do not share
vertices). A matching is called {\em perfect} if its cardinality
is maximal.
\end{defn}

\smallskip
\noindent The following proposition shows the existence of a
unique perfect matching of $H_D$.

\begin{thm} \label{HDbip}
Let $D$ be a digraph. The undirected graph $H_D$ is bipartite and
it has a perfect matching.
\end{thm}
\textbf{Proof}. Let $V=\{v_1, \ldots, v_n\}$ and $Z=\{z_1, \ldots,
z_n\}$. $V$ and $Z$ are disjoint subsets of $V(H_D)$, whose union
is exactly $V(H_D)$. By definition of $E(H_D)$ every edge in $H_D$
connects a vertex in $V$ with a vertex in $Z$. So $H_D$ is
bipartite with bipartition sets $V$ and $Z$, by definition of
bipartite graph. Finally the edges $f_i=\{z_i,v_i\}$ for $i=1,
\ldots, n$ are independent and they are $n$ in a graph with $2n$
vertices. So the set $M=\{f_1, \ldots, f_n\}$ is a perfect
matching for $H_D$ by its own definition. \qed

\smallskip
\noindent
 Now, given a bipartite graph $G$, we want to find a digraph $D$, with $H_D=G$.
 The following theorem shows some results in this direction.

\begin{thm}
Let $G$ be a bipartite undirected graph with $V(G)=\{v_1,\ldots,
v_n,\\ z_1,\ldots,z_n\}$. Then the following propositions are
equivalent:
\begin{enumerate}
\item $G$ has a perfect matching; \item there exists a digraph
$D$, with $V(D)=\{v_1, \ldots, v_n\}$ such that $H_D=G$.
\end{enumerate}
\end{thm}
\textbf{Proof}. (1) $\Rightarrow$ (2). If $G$ has a perfect
matching, then it has $n$ edges, that share no vertices. We can
relabel the vertices in such a way as the $n$ independent edges
are $f_i=\{z_i, v_i\}$, for $i=1, \ldots, n$. So $V=\{v_1, \ldots,
v_n\}$ and $Z=\{z_1, \ldots, z_n\}$ are the two bipartition sets.
Otherwise either a vertex $z_i$ should be in a bipartition set
with some $v_j$, (or a vertex $v_i$ should be in a bipartition set
with some $z_j$). This fact would imply that the edge $f_i$ has
two vertices in the same bipartition set, against the definition
of bipartite graph. Now it is sufficient to define $D$ in such a
way as $V(D)=\{v_1, \ldots, v_n\}$ and $[v_i,v_j]\in
E(D)$ whenever $\{z_i, v_j\} \in E(G)$.\\
(2) $\Rightarrow$ (1). It follows from theorem \ref{HDbip}. \qed

\section{Cycles in $D$ and $H_{D}$}
Let $G$ be an undirected graph. There exists a binomial ideal
$I(G,E)$ in the vertices and edges associated to $G$ (\cite{CF04}),
while there exists a binomial ideal $I(E)_{E} \subseteq I(G,E)$ in
the edges associated to the even closed walks in $G$ (\cite {HO98}).
Given an even closed walk
$C=(e_{i_{1}}=\{v_{i_{1}},v_{i_{2}}\},\ldots,
e_{i_{2q-1}}=\{v_{i_{2q-1}},v_{i_{2q}}\},
e_{i_{2q}}=\\\{v_{i_{2q}},v_{i_{1}}\})$ of $G$, let
$f_{C}=\prod_{k=1,\ldots,q}e_{i_{2k-1}}-\prod_{k=1,\ldots,q}e_{i_{2k}}$
be the corresponding binomial in $I(E)_{G}$.\\
The following theorem  is a well known result about even closed
walks in an undirected graph.

\begin{thm}\label{Hibi}$($\cite{Vi95}$)$,$($\cite{HO98}$)$\\
If $G$ is an undirected graph, then the toric ideal $I(E)_{G}$,
associated to its incidence matrix, is generated by all  binomials
$f_C$, where $C$ is an even closed walk of $G$.
\end{thm}

\smallskip
\noindent Now we want to extend this result to the case of
digraphs.

\begin{thm}\label{cycles}
Let $D$ be a digraph. The toric ideal $I(E)_D$ is generated by all
binomials $f_{C}$, where $C$ is a cycle of $D$.
\end{thm}
\textbf{Proof}. By looking at the undirected  graph $H_D$
associated to $D$, then the ideal $I(E)_{H_D}$ = $I(H_D, E) \cap
K[e_{1},\ldots,e_{m},f_{1},\ldots,f_{n}]$ is generated by all even
closed walks in $H_D$ and then by all even cycles, because $H_D$
is bipartite by theorem \ref{HDbip}. Now $I(D,E)=\pi(I(H_D,E))$
and the binomial edge ideal $I(E)_D$=$I(D,E)$ $\cap$
$K[e_{1},\ldots,\\e_{m}]$ is equal to $\pi(I(H_D,E)) \cap
K[e_{1},\ldots,e_{m},f_{1},\ldots,f_{n}]$. It follows that
$I(E)_D$ is generated by  binomials $f$=$\prod_{i\in
I}f_{i}\prod_{i'\in I'}e_{i'}-\prod_{j\in J}f_{j}\prod_{j'\in
J'}e_{j'}$  with $I, J \subseteq \{1,\ldots,n\}$, $I', J'
\subseteq \{1,\ldots,m\}$, $|I|+|I'|=|J|+|J'|$ and
$f_{i}$=$f_{j}$=$1$ for all $i,j=1,\ldots,n$. If $C$ is a direct
cycle in $D$, then we can relabel the vertices and the edges is
such a way as  $C$= \{$e_{1}=[v_{1},v_{2}]$,
\ldots,$e_{q-1}=[v_{q-1},v_{q}]$,$e_{q}=[v_{q},v_{1}]\}$ and
$C'$=$\{e_{1},f_{2},e_{2},f_{3}, \ldots, e_{q-1},f_{q},
e_{q},f_{1}\}$ is an even cycle in $H_D$. So the binomial
$f_{C'}$=$\prod_{h=1,\ldots,q}e_{h}- \prod_{h=1,\ldots,q}f_{h}$ is
in  $I(H_D,E) \cap K[e_{1}, \ldots, e_{m}, f_{1}, \ldots, f_{n}]$
and the binomial $f_{C}$ = $\prod_{h=1,\ldots,q}e_{h}-1$ is in
$I(D,E) \cap
K[e_{1},\ldots,e_{m}]$.\\
Now let $C$ be an undirected cycle in $D$. Once again we can
relabel the vertices and the edges in such a way as $C$=
$\{\{v_{1},v_{2}\}, \ldots, \{v_{q-1},
v_{q}\},\{v_{q+1}=v_{1}\}\}$. Let $C'$ be the path in $H_D$ given
in the following way: $C'$=$g_{1},\ldots,g_{q}$, where
$g_{h}=e_{h}f_{h+1}$ if $e_{h}=\{z_{h},v_{h+1}\}$ and
$[v_{h},v_{h+1}] \in E(D)$, while $g_{h}=f_{h+1}e_{h}f_{h}$ if
$e_{h}=\{z_{h+1},v_{h}\}$ and $[v_{h+1},v_{h}] \in E(D)$.
$C'$=$\{z_{1},v_{2},z_{2},v_{3},
\ldots,z_{q-1},v_{q},z_{q},v_{1}\}$ is an even cycle in $H_D$. The
corresponding binomial $f_{C'}$ is in the ideal $I(H_D,E) \cap
K[e_{1},\ldots,e_{m},f_{1},\ldots,f_{n}]$, while its image in
$K[e_{1},\ldots,e_{m}]$ is the binomial $f_{C}$
corresponding to $C$.\\
Conversely given a cycle $C'$ in $H_D$, then it is an even cycle,
because $H_D$ is bipartite. So
$C'$=$\{z_{i(1)},v_{i(2)},z_{i(3)},v_{i(4)},
\ldots,z_{i(q-3)},v_{i(q-2)},z_{i(q-1)},v_{i(q)},z_{i(q+1)}$=$z_{i(1)}\}$.
Let $e_{h}=\{z_{i(h)},v_{i(h+1)}\}$ whenever $i(h) \neq i(h+1)$
and let $f_{h}=\{z_{i(h)},v_{i(h+1)}\}$ whenever $i(h)=i(h+1)$.
$C'$ determines the cycle $C$= $\{v_{i(1)},v_{i(2)}, \ldots,
v_{i(q)},v_{i(q+1)}=v_{i(1)}\}$, where $v_{i(h)}=v_{i(h+1)}$,
whenever $i(h)=i(h+1)$ and $e_{i(h)}=\{v_{i(h)}, v_{i(h+1)}\}$ is
in $E(D)$, whenever $i(h) \neq i(h+1)$. \qed

\begin{rem}
The existence of the toric ideal as above was proved in different
way by using circuits associated to the transpose of the incidence
matrix in 2005 by Reyes (\cite{Re05}) and in 2006 by Gitler, Reyes
and Villarreal (\cite{GRV06}.
\end{rem}
\begin{cor}\label{cycles undirected}
Let $G$ be an undirected graph and let $G^d$ be a directed graph,
whose underlying graph is $G$ (e.g. let $G^d$ be a directed graph
such that $G^d_{u}=G$). Then the ideal $I(G^d, E)$ is generated by
all binomials $f_C$, such that $C$ is a cycle of $G$.
\end{cor}
\textbf{Proof}. Let $V(G)=\{v_1, \ldots, v_n\}$ and
$E(G)=\{e_{i,j}=\{v_i,v_j\}, i,j \in \{1, \ldots, n\}\}$. It is
possible to construct such directed graph $G^d$ by setting
randomly a direction in each edge. Let  $V(G^d)=V(G)$ and let
$E(G^d)=\{e_{i,j}=[v_i,v_j]$, such that $\{v_i,v_j\} \in E(G)\}$.
By theorem as above $I(G^d, E)$ is generated by binomial
$f_C=\prod_{i \in I} e_i - \prod_{j \in J} e_j$, where $C$ is a
cycle of $G^d$ and $J=\emptyset$ whenever $C$ is a directed cycle.
Since every cycle in $G^d$  is a cycle in $G$, then we have the
thesis. \qed

\noindent The theorem as above gives also decision procedures for
the existence of directed and undirected cycles in a digraph as in
the following remark.

\begin{rem}\label{diundicycles}
The  proof of the theorem as above allows to check whether a given
digraph $D$ has cycles  and  whether the cycles are either directed
or undirected. By using Maple 10 we implemented the corresponding
decision procedures. The  algorithm works as follows: given a
digraph $D$ by using the package {\em networks} we find the
incidence matrix $M$ of $D$ and then by using the package {\em
Groebner} we get a Gr\"{o}bner basis $B$ of the toric ideal $I(E)_D$
of $M$. If $B$ contains a polynomial of the form $f=\prod_{i \in I}
e_i -1$, then in $D$ there is a directed cycle of length $|I|$.
Furthermore if $B$ contains a polynomial of the form $g=\prod_{i \in
I} e_i - \prod_{j \in J} e_j$, then in $D$ there is an undirected
cycle of length $|I|+|J|$. Of course the binomials associated to
even and odd cycles in a digraph are in the edge toric ideal, while
in the undirected case only even cycles appear in the edge toric
ideal associated to the graph. By using corollary 1 we are able to
check whether an undirected graph $G$ has even and odd cycles. We
can construct $G^{d}$ by simply choosing randomly a direction for
each edge and we can find the toric ideal $I(G^{d},E)$. Such ideal
is generated by binomials in the form either $f=\prod_{i \in
I}e_{i}-1$ or $g=\prod_{i \in I}e_{i}-\prod_{j \in J}e_{j}$, that
are in one to one correspondence with cycles in $G$.
\end{rem}

\noindent By using the algorithm sketched in the last remark it is
possible to deduce properties about the digraphs $D_1$ and $D_2$ in
examples \ref{D1} and \ref{D2}.

\begin{ex}\label{tD1}
The toric ideal basis of the ideal $I(E)_{D_1}$, with respect to
the lexicographic term order $\sigma_{1}$ with
$v_{1}>_{\sigma_{1}}
z_{1}>_{\sigma_{1}}v_{2}>_{\sigma_{1}}z_{2}>_{\sigma_{1}}v_{3}>_{\sigma_{1}}z_{3}
>_{\sigma_{1}}v_{4}>_{\sigma_{1}}z_{4}>_{\sigma_{1}}v_{5}>_{\sigma_{1}}
z_{5}>_{\sigma_{1}}e_{3}>_{\sigma_{1}}e_{1}>_{\sigma_{1}}e_{4}
>_{\sigma_{1}}e_{2}>_{\sigma_{1}}e_{5}>_{\sigma_{1}}e_6$ is
$(e_1 e_2 e_3-1, e_1 e_2 e_5 - e_4, e_3 e_4-e_5)$. So  the digraph
$D_1$ has a directed cycle with the edges $e_1, e_2, e_3$ and two
undirected cycles with the edges $e_3, e_4, e_5$ and $e_1, e_2,
e_4, e_5$.
\end{ex}

\begin{ex}\label{tD2}
The digraph $D_2$ is a DAG, in fact the toric ideal basis of the
ideal $I(E)_{D_2}$, with respect to the lexicographic term order
$\sigma_{2}$ with $v_{1}>_{\sigma_{2}}
z_{1}>_{\sigma_{2}}v_{2}>_{\sigma_{2}}z_{2}>_{\sigma_{2}}v_{3}>_{\sigma_{2}}z_{3}
>_{\sigma_{2}}v_{4}>_{\sigma_{2}}z_{4}>_{\sigma_{2}}v_{5}>_{\sigma_{2}}
z_{5}>_{\sigma_{2}}e_{1}>_{\sigma_{2}}e_{3}>_{\sigma_{2}}e_{2}
>_{\sigma_{2}}e_{4}>_{\sigma_{2}}e_{5}$ is
$(e_1 e_2 - e_3)$ and the only cycle with the edges $e_1, e_2,
e_3$ is undirected.
\end{ex}

\smallskip
\noindent
Now we study a property of digraphs, that can be easily
checked by using theorem \ref{cycles}.

\begin{defn}
Given a digraph $D$, a vertex $v$ is reachable from another vertex
$u$ if there is a directed path that starts from $u$ and ends at
$v$. If $v$ is reachable from $u$, then $u$ is a {\em predecessor}
of $v$ and $v$ is a {\em successor} of $u$.
\end{defn}

\begin{defn}
Let $D$ be a digraph and let $v$ and $u$ be vertices in $D$. $D$
is  a $UPD$ (Unique Path Digraph) iff whenever $v$ is a successor
of $u$, then there is a unique elementary path between $u$ and
$v$.
\end{defn}

\smallskip
\noindent It is easy to show that a cycle of a $UPD$ is directed.
The following theorem is useful for a characterization of a $UPD$.

\begin{thm}\label{UPD}
Let $D$ be a digraph and let $C_1$ and $C_2$ be two directed
cycles in $D$, such that $C_1$ and $C_2$ are not edge-disjoint.
%$C_1 \cap C_2 \neq \emptyset$.
%\footnote{The intersection of two cycles is the set of edges in
%both  cycles.}
Then $C=(C_1 \cup C_2) \backslash (C_1 \cap C_2)$
is an undirected cycle.
\end{thm}

\textbf{Proof}. Let $I(E)_D$ be the edge ideal associated to $D$.
By theorem \ref{cycles} the  binomials representing the two
directed cycles $C_1$ and $C_2$ lie in $I(E)_D$. Let $f_1=\prod_{i
\in I} e_i -1$ be the binomial associated to $C_1$ and let
$f_2=\prod_{j \in J} e_j -1$ be the binomial associated to $C_2$.
Since $C_1 \cap C_2 \neq \emptyset$, then $K=I \cap J \neq
\emptyset$ and  $f_1=(\prod_{k \in K} e_k \cdot \prod_{i \in I
\backslash K} e_i) -1$ while $f_2=(\prod_{k \in K} e_k \cdot
\prod_{j \in J \backslash K} e_j) -1$. Let $\sigma$ be a
lexicographic term ordering in $\mathbb{K}[e_1, \ldots, e_m]$,
where $m$ is the number of edges in $D$. The S-polynomial between
$f_1$ and $f_2$ with respect to the term ordering $\sigma$ is the
binomial $f= -\prod_{j \in J \backslash K} e_j + \prod_{i \in I
\backslash K} e_i$. By remark \ref{diundicycles}  the cycle $C$
corresponding to $f$ is an undirected cycle lying in $D$ and
$C$=$(C_1 \cup C_2) \backslash (C_1 \cap C_2)$. \qed

\begin{cor}\ref{UPD}
Let $D$ be a digraph. $D$ is a $UPD$ if and only if its  cycles
are directed  with at most a vertex in common.
\end{cor}

\textbf{Proof}. If $D$ is a tree or a forest there is nothing to
prove. Otherwise all  cycles have to be directed and by theorem
\ref{UPD} they can have at most a vertex in common. \qed

\begin{rem}
The UPD property of a digraph $D$ can be easily checked by using
theorem \ref{cycles} and corollary \ref{UPD}. In fact it is easy
to show that $D$ is $UPD$ if and only if either the binomial edge
ideal $I(E)_{D}=(0)$, e.g. $D$ has no cycles, or $I(E)_{D}$ is
generated by binomials $f_{h}=\prod_{j=1,\ldots,k(h)}e_{j}-1$ with
$\prod_{j=1,\ldots,k(h)}e_{j}$ and $\prod_{j=1,\ldots,k(h')}e_{j}$
coprime for all $h \neq h'$.
\end{rem}

\section{Linear ideals arising from digraphs}
\label{linear}

Here we will introduce some linear ideals, that we can  associate to
a digraph and we can use in decision procedures. By using algorithms
from linear algebra, we could loose  some property of a digraph.
Anyway such algorithms are very useful, because they are fast.

\noindent Let $D$ be a digraph and let $H_D$ be the associated
undirected graph defined as before. In \cite{CF04} the extended
linear edge ideal (respectively the linear edge ideal) can be
associated to $H_D$ and relations between the extended linear edge
ideal and the extended binomial edge ideal (respectively the linear
edge ideal and the binomial edge ideal) are shown.

\smallskip
\noindent Let $R=K[v_{1},z_{1},f_{1}, \ldots, v_{n},z_{n},f_{n},
e_{1}, \ldots, e_{m}]$ and let $\pi': R \rightarrow
R/(f_{1},\ldots f_{n})$ be the canonical ring homomorphism defined
by $\pi'(v_{i})=v_{i}$, $\pi'(z_{i})=z_{i}$, $\pi'(f_{i})=0$,
$\pi'(e_{j})=e_{j}$ for all $i=1,\ldots,n$ and $j=1,\ldots,m$.

\begin{defn}
The ideal $LI(D,E)$=$\pi'(LI(H_D,E))$ in $K[e_1, \ldots, e_m, v_1,
\ldots, v_n]$ is the  {\em linear extended edge ideal of} $D$. The
ideal $LI(E)_D$ = $LI(D,E)\cap K[e_1, \ldots, e_m]$ is the {\em
linear edge ideal of} $D$.
\end{defn}

\begin{rem}
It is easy to show that the ideal $LI(D,E)$=$($ $e_h+v_{i}-v_{j}$:
$e_h=[v_{i},v_{j}]$ in E(D)$)$ in $K[e_1, \ldots, e_m, v_1,
\ldots, v_n]$ by definition of $\pi'$.
\end{rem}

\noindent By repeating the proof as in \cite{CF04} if we take the
matrix $M=IM(D)$ and the linear homomorphism of semigroup algebras
$\psi$ : $K[e_1,\ldots,e_m] \longrightarrow
K[v_1,\ldots,v_n,v_{1}^{-1},\ldots,\\v_{m}^{-1}]$ defined by
$\psi(e_h)=v_j-v_i$ whenever $e_h=[v_i,v_j]$ for all
$h=1,\ldots,m$, then $LI(D,E)$ coincides with the ideal
$LI_{M}$=$ker(\psi)$  in $K[e_1,\ldots,e_m]$, that is called the
{\em linear ideal of} $M$.

\begin{ex}
Let $D_1$ be the graph as in example \ref{D1} and \ref{tD1}. The
linear extended edge ideal of $D_1$ is $LI(D_1,E)=(e_1+v_1-v_2,
e_2+v_2-v_3, e_3+v_3-v_1, e_4+v_1-v_4, e_5+v_3-v_4, e_6+v_3-v_5)$,
while the linear edge ideal of $D_1$, is
$LI(E)_{D_1}=(e_3+e_4-e_5,e_1+e_5+e_2-e_4)$.
\end{ex}

\begin{ex}
Let $D_2$ be the graph as in example \ref{D2} and \ref{tD2}. The
linear extended edge ideal of $D_2$ is $LI(D_2,E)=(e_1+v_1-v_2,
e_2+v_2-v_3, e_3+v_1-v_3, e_4+v_4-v_3, e_5+v_3-v_5)$, while the
linear edge ideal of $D_2$, is $LI(E)_{D_2}=(e_1-e_3+e_2)$.
\end{ex}

\smallskip \noindent In order to show the relations between the
linear ideal and the toric ideal associated to the incidence
matrix $IM(D)$ we need the following definition.

\begin{defn}
{\em (\cite{St95})} Let $M$=$(m_{ij})_{i=1,\ldots,m,j=1,\ldots,n}$
be a $(m,n)$-matrix with $m_{ij}$ in ${\bf N}_{0}$ for all $i,j$
and let $\pi_{M}$ : ${\bf N}_{0}^{n} \longrightarrow {\bf Z}^{m}$
and $\pi_{M}'$ : ${\bf Z}^{n} \longrightarrow {\bf Z}^{m}$ be the
semigroup homomorphisms as in definitions \ref{deftoric}. $Ker(M)$
is the kernel of the semigroup homomorphism $\pi_{M}'$. Moreover
if a finite set $F$ generates $Ker(M)$ as {\bf Z}-module, then
$I(F)$= $(e^{u^{+}}-e^{u^{-}}, u \in F)$ is a {\em lattice ideal
associated to} $M$.
\end{defn}

\begin{rem}\label{rem5.1}
 $\phi_{n}^{-1}(Ker(M))$ is a {\bf Z}-submodule of
$\bigoplus_{j=1,\ldots,n}{\bf Z}X_{j}$ and it is the kernel of the
{\bf Z}-module homomorphism $\phi_{m}\pi_{M}'\phi_{n}$ by definition
of the homomorphisms $\pi_{M}'$, $\phi_{n}$ and $\phi_{m}$. Finally
since $\bigoplus_{j=1,\ldots,n}{\bf Z}X_{j}$ is a {\bf Z}-submodule
of $A=K[X_1,\ldots,X_n]$, then $\phi_{n}^{-1}(Ker(M))$=$ker(\psi)
\cap \bigoplus_{j=1,\ldots,n}{\bf Z}X_{j}$.
\end{rem}

\smallskip
\noindent The following definition is also useful and
related to the notion of saturation.

\begin{defn}
Let $I$ be an ideal in the polynomial ring $A$ and let $f \in A$.
$I:f^{\infty}$=($g$: $gf^{m} \in I$ for some $m \in {\bf N}_{0}$).
\end{defn}

\smallskip
\noindent Of course $I \subseteq I:f^{\infty}$ for all $f \in A$

\smallskip
\noindent The relation between the ideals $LI(E)_D$ and
$I(E)_D$ is
given by the following facts.\\
First of all in  \cite{HS95} and (\cite{St95}, p.114) it is shown
that if a finite set $B$ generates $Ker(IM(D))$ as {\bf Z}-module,
$J_{0}$=$I(B)$=$(e^{u^{+}}-e^{u^{-}}, u \in ker(IM(D))$ and
$J_{i}$=$(J_{i-1}:e_{i}^\infty)$ for all $i=1,\ldots,m$, then
$J_{m}=I(D)_{E(D)}$. \\
Furthermore the generators of $J_{0}$ are in one to one
correspondence with the elements of $F$ and then in one to one
correspondence with the elements of $\phi_{n}^{-1}(F)$, that
generate the {\bf Z}-module $\phi_{n}^{-1}(Ker(IM(D)))$=$LI(IM(D))
\cap \bigoplus_{j=1,\ldots,n}{\bf Z}e_{j}$=$LI(D)_{E} \cap
\bigoplus_{j=1,\ldots,n}{\bf Z}e_{j}$ by remark \ref{rem5.1}.

\begin{rem}
The algorithm in \cite {HS95} called the {\em saturation
algorithm} is   one of the existing algorithms for finding such
ideal. Another well known algorithm is the {\em Lift-and-Project
Algorithm} in \cite{BLR99} and recently the algorithm in
\cite{HM06}.
\end{rem}

\smallskip
\noindent Now  the relation between the ideals $LI(D,E)$ and
$I(D,E)$ is the same as the relation between the ideals $LI(G,E)$
and $I(G,E)$ when $G$ is an undirected graph as in \cite{CF04}.
The  relations between binomial and linear ideals associated to a
digraph give also the following theorem.

\begin{thm}\label{licycles}
Let $D$=$(V(D),E(D))$ be a simple digraph.\\
(i) The directed cycle
$C$=$(e_{1}, \ldots, e_{k})$ is in $D$ iff the linear
polynomial $h_{C}$=$\sum_{h=1, \ldots, k}e_{h}$ is in $LI(D,E)$.\\
(ii) The undirected cycle $C$=$(e_1, \ldots, e_k)$ is in $D$ iff
the linear polynomial $h_{C}$=$\sum_{i \in I} e_i - \sum_{j \in J}
e_j$, with $|I|+|J|=k$ is in $LI(D)_{E}$.
%(iii)Let $\sigma$ be a lexicographic term ordering on the set of the
%power products in $\{e_1,\ldots,e_m,v_1,\ldots,v_n\}$ with $v_i >
%e_j$ for all $i$ and $j$ and
%$v_{2q}>_{\sigma}v_{2q-1}>_{\sigma}\ldots>_{\sigma}v_{1}$. If $C$ is
%minimal, then the polynomial $h_{C}$ is in a Gr\"{o}bner basis of
%$LI(G,E)$ with respect to $\sigma$.
\end{thm}
\textbf{Proof}. The proof of (i) and (ii)
%and ?(iii)
is the same as the proof of theorems about the corresponding
binomial ideals.

\begin{rem}
The theorem as above shows that the hypothesis on the
characteristic of the field $K$ is necessary. In fact the theorem
does not hold when the characteristic of $K$ is different from
$0$.
\end{rem}

\begin{ex}\label{lD1}
Let $D_1$ be the graph as in examples \ref{D1} and \ref{tD1}. Then
the Gr\"{o}bner basis of $LI(D_1,E)$ with respect to $\sigma_1$ is
$(e_4-e_5+e_3, -e_4+e_5+e_1+e_2, -e_5+v_4+e_6-v_5, e_6+v_3-v_5,
v_2+e_2-v_5+e_6, v_1-v_5+e_4+e_6-e_5)$.
\end{ex}

\begin{ex}\label{lD2}
Let $D_2$ be the graph as in examples \ref{D2} and \ref{tD2}. The
Gr\"obner basis of the ideal $LI(D_2,E)$ with respect to
$\sigma_2$ is $(e_1-e_3+e_2,e_5-v_5+e_4+v_4, e_3+v_1+e_5-v_5,
e_5-v_5+v_3, v_2- v_5+e_2+e_5)$.
\end{ex}

\noindent By using linear ideals, the procedures as above are very
fast. Now  theorem \ref{licycles} establishes only that a
polynomial $p$, that corresponds to a directed cycle, is in the
ideal. There is no proof that $p$ is in some Gr\"obner basis of
$LI(E)_D$ as in example \ref{lD1} with $p=e_1+e_2+e_3$ and then we
cannot decide if $D$ is a DAG. Anyway we can use the linear
ideals, once we want to check the existence of cycles.

\medskip
\noindent The same results can be obtained by using the classical
tools of computational graph theory applied to cycle spaces and
cycle bases as in  \cite{KMMP04} and  \cite{KM05}.

\begin{defn}
Let $D$=$(V,E)$ be a simple digraph and let $C$ be  cycle in $D$.
The {\em incidence vector of} $C$ is a vector in
$\{-1,0,1\}^{|E|}$ defined as follow. For every $e \in E$ C(e) is
equal to $1$ if $e=[v_i,v_j] \in C$ , C(e) is equal to $-1$ if
$-e=[v_j,v_i] \in C$  and C(e) is equal to $0$ if $e \notin C$.
The {\em cycle space of} $D$ is the vector space over ${\bf Q}$
spanned by the incidence vectors of its cycles. A {\em cycle basis
of} $D$ is a basis of the cycle space.
\end{defn}

\begin{rem}
If $D$ is connected, then the cycle space of $D$ has dimension
$d=|E|-|A|+1$.
\end{rem}

It is well known (\cite{GLS03} and \cite{Bo98}) that if $D$ is a
connected digraph, then the cycle space of $D$ is the kernel of
the linear map from the edge space $E(D)$ to the vertices space
$V(D)$ defined by the incidence matrix $IM(D)$ and a basis of such
vector space is a cycle basis of $D$. Furthermore if $D$ is
strongly connected, i.e. for every pair of vertices $u$ and $v$ in
$D$, there are directed paths both from $u$ to $v$ and from $v$ to
$u$, then it is possible to find a cycle basis that contains all
cycles of length $2$ and all directed cycles.

\medskip
\noindent Now let $D=V(D),E(D)$ be a simple digraph. Let
$n=|V(D)|$ and let $m=|E(D)|$. By looking at the ideal $LI(D,E)$,
then it is generated by the kernel of the linear map associated to
the matrix $(m,n+m)$-matrix $(IM(D)^t,-I(m))$, while $LI(E)_D$ is
generated by a cycle basis of $D$, when we identify the incidence
vector of a cycle $C$ with the corresponding linear polynomial in
the edges in $C$. The results in \cite{GLS03} shows that we can
know all directed cycles in a digraph when some other hypothesis
are satisfied, (in particular digraphs with double edges, i.e.
connecting two vertices in both directions, are allowed) while a
basis of the toric ideal associated to the transpose of the
incidence matrix, which is also an ideal in the edges containing
the lattice ideal associated to the same matrix, contains a
minimal set of directed cycles in $D$. In fact in general the set
of all directed cycles in a digraph $D$ is not a basis of a vector
space, because the sum of two directed cycles is not well defined.

\smallskip
\noindent
 There are many digraphs that are not strongly
connected and we cannot apply neither the algorithm in
\cite{GLS03} nor our theorem \ref{licycles}. So we have to use
theorem \ref{cycles}. Finally in \cite{KMMP04} the authors give a
fast algorithm to compute a cycle basis of minimal weight in
undirected graphs.

\section{Sink and source covers of a digraph}
Here we introduce another undirected graph $K_{D}$, that we can
associate to a digraph $D$. We can prove some properties of $D$
trough the properties of such graph.

\begin{defn}
Let $D$ be a digraph. Let $K_D$ be the undirected subgraph of
$H_D$ with $V(K_D)$=$V(D) \cup \{z_{1},\ldots,z_{n}\}$=
 and $E(K_D)$=\{$e=\{z_{i},v_{j}\}$: $[v_{i},v_{j}]
\in E(D)$\}. $K_D$ is called the {\em sink-source undirected graph} associated to $D$.
\end{defn}

\smallskip
\noindent In  \cite{CF04} it is introduced the extended vertex
ideal of the undirected graph $G_u$  $I(G_u,V(G_u))$=($v_i-\prod
e_{h}$: $v_i$ is in  $e_h \in E(G_u)$).

\begin{defn}
Let $D$ be a digraph and let $D^*=K_D \setminus L$, where $L$ is the
set of isolated vertices. The {\em extended divertex} ideal of $D$
is equal to $I(D^*,V)$.
\end{defn}

\smallskip
\noindent The following result  can be found in \cite{CF04}).

\begin{thm}\label{bip undirected}
%\label{bip undirected} (Carr\`a Ferro and Ferrarello, 2004.
%See \cite{CF04})\\
Let $G$=$(V(G),E(G))$ be a simple undirected connected graph without
isolated vertices and let $I(G,V)$ be the extended vertex ideal of
$G$. Then the ideal $I(V)_{G}$=$I(G,V)\cap K[v_1,\ldots,v_n]$
contains an irreducible polynomial $p$ of the form $p$=$\prod_{j \in
J}v_{j}-\prod_{k \in K}v_k$, if and only if $G$ is bipartite.
Moreover the partition sets are $V'$=$\{v_j: j \in J\}$ and
$V''$=$\{v_k: k \in K\}$.
\end{thm}

\smallskip
\noindent A possible generalization of the last theorem to the
digraph uses the notion of directly bipartite graph.

\begin{thm}
Let $D$=$(V(D),E(D))$ be a  digraph and let $I(V)_{D}$ be the
divertex ideal of $D$. Then $D$ is directly bipartite if and only if
$I(V)_{D}$ contains an irreducible polynomial of the form
$p=\prod_{i \in I} z_i - \prod_{j \in J} v_j$, with $I \cap
J=\emptyset$, $I \cup J= \{1, \ldots, n\}$.
\end{thm}
\textbf{Proof}. $D$ is directly bipartite if and only if $V(D)=V_I
\cup V_J$, with $V_I=\{v_i: \\i \in I\}$, $V_J=\{v_j: j \in J\}$,
$I \cap J=\emptyset$, $I \cup J=\{1, \ldots, n\}$ and every edge
in $E(D)$ goes from a vertex in $V_I$ to a vertex in $V_J$. By
definition of $K_D$ this last fact is equivalent to say that $K_D$
is bipartite with partition sets $Z=\{z_i: i \in I\}$ and
$V=\{v_j: j \in J\}$ and with the $n$ isolated vertices $\{z_j: j
\in J\} \cup \{v_i: i \in I\}$. Now by theorem \ref{bip
undirected} this fact is equivalent to say that the binomial
$p=\prod_{i \in I} z_i - \prod_{j \in J} v_j$ is in the extended
vertex ideal of $K_D$. Finally, this last ideal coincides with the
extended vertex ideal of $D$ by definition. \qed

\smallskip
\noindent Let us recall that a {\em vertex cover} $W$ of an
undirected graph $G=(V,E)$ is a subset of $V$, such that every
edge in $E$ is incident with at least a vertex in $W$. A vertex
cover $W$ is called {\em minimal} if every subset of $W$ is not a
vertex cover. It is possible to generalize the concept of minimal
vertex cover for digraphs, in the following way:

\begin{defn}
Let $G=(V,E)$ be a digraph. A {\em source cover of} $G$ is a vertex
cover $V'$ of $G$, such that every edge in $G$ leaves every vertex
in $V'$. A source cover $V'$ of $G$ is called {\em minimal} if no
subset of
$V'$ is a source cover of $G$.\\
A {\em sink cover of} $G$ is a vertex cover $V'$ of $G$, such that
every edge in $G$ does not leave any vertex in $V'$. A sink cover
$V'$ of $G$ is called {\em minimal} if no subset of $V'$ is a sink
cover of $G$.
\end{defn}

\noindent The following proposition is the key result for our
purposes.

\begin{prop} \label{Villarreal} $($ \cite{Vi01}$)$
 Let $K[v]=K[v_1, ..., v_n]$ be a polynomial ring over a field K
and let $G$ be an undirected graph. If {\em P} is the ideal of
$K[v]$ generated by $A=\{v_{i1}, \ldots, v_{ir}\}$, then {\em P}
is a minimal prime ideal containing  the edge ideal $I(G)_{E}$ if
and only if $A$ is a minimal vertex cover of $G$.
\end{prop}

\smallskip
\noindent By using the undirected graph $K_G$ it is possible to
find source and sink covers of a digraph $G$, according to the
following proposition.

\begin{prop}\label{sscovers}
Let $G$ be a directed graph and let $K_G$ be the associated
source-sink undirected graph.
Let $V'$ be a vertex cover of $K_G$. $V'=\{v_{i(1)},\ldots,v_{i(l)}\}$ is a source cover of $G$ if and %%@
only if $\{z_{i(1)},\ldots,z_{i(l)}\}$ is a vertex cover of $K_G$.
In similar way
$V'=\{v_{i(1)},\ldots,v_{i(l)}\}$ is a sink cover of $G$ if and only if $\{v_{i(1)},\ldots,v_{i(l)}\}$ %%@
is a vertex cover of $K_G$.
\end{prop}
\textbf{Proof}. Let $V'=\{v_{i(1)},\ldots,v_{i(l)}\}$ be a source
cover of $G$.
So $V'$ is a vertex cover of $G$, such that each edge of $G$ leaves at least one vertex in $V'$. It %%@
follows that $\{z_{i(1)},\ldots,z_{i(l)}\}$ is a vertex cover of $K_G$ by its own definition. In %%@
similar way if $V'$ is a sink cover of $G$, then
$V'$ is a vertex cover of $G$, such that each edge of $G$ does not leave any   vertex in $V'$. So $V'$  %%@
is a vertex cover of $K_G$ by its own definition. Conversely, let $Z'$=$\{z_{i(1)},\ldots,z_{i(l)}\}$ %%@
be a vertex cover of $K_G$. So each edge in $K_G$ is incident with at least one vertex in $Z'$. It %%@
follows that each edge in $G$ is incident with at least one vertex in  %%@
$V'=\{v_{i(1)},\ldots,v_{i(l)}\}$  and it leaves such vertex by definition of $K_G$. So $V'$ is a %%@
source cover of $G$. In similar way if $V'=\{v_{i(1)},\ldots,v_{i(l)}\}$ is a vertex cover of $K_G$, %%@
then it is a vertex cover of $G$. Moreover $V'$ is a sink cover of
$G$ by definition of $K_G$. \qed

\begin{ex}
Let $D_1$ be the graph as in example \ref{D1}. $\{z_1,z_2,z_3\}$
and $\{v_1,v_2,v_3,v_4,\\v_5\}$ are minimal vertex covers of
$K_{D_1}$. $\{v_1,v_2,v_3\}$ is a source cover of $D_1$ and
$\{v_1,v_2,v_3,v_4,v_5\}$ is a sink cover of $D_1$. $D_1$ is not
directly bipartite.
\end{ex}

\begin{ex}
Let $D_2$ be the graph as in example \ref{D2}. $\{z_1,z_2,z_3,z_4\}$
and $\{v_2,v_3,v_5\}$ are minimal vertex covers of $K_{D_2}$.
$\{v_1,v_2,v_3,v_4\}$ is a source cover of $D_2$ and
$\{v_2,v_3,v_5\}$ is a sink cover of $D_2$. $D_2$ is not directly
bipartite.
\end{ex}

%\begin{ex}
%Let $D_3$ be the graph as in example \ref{D3}. $\{z_1,z_3,z_4\}$ and
%$\{v_1$, $v_2$, $v_3$, $v_4\}$ are minimal vertex covers of
%$K_{D_3}$. $\{v_1,v_3,v_4\}$ is a source cover of $D_3$ and
%$\{v_1,v_2,v_3,v_4\}$ is a sink cover of $D_3$. $D_3$ is not
%directly bipartite.
%\end{ex}

\begin{ex}
Let $D_4$ be the digraph with vertex set $V(D_4)=\{v_1, \ldots,
v_5\}$ and edge set $E(D_4)$ = $\{e_1=[v_1,v_2]$, $e_2=[v_3,v_2]$,
$e_3=[v_1,v_4]$, $e_4=[v_3, v_4]$, $e_5=[v_3, v_5]$. $D_4$ is
directly bipartite and it has a sink cover and a source cover. In
fact the divertex ideal of $D_4$, obtained by intersecting the
ideal $(z_1-e_1 e_3, v_2-e_1 e_2, z_3-e_2 e_4 e_5, v_4-e_3 e_4,
v_5-e_5)$ with $K[z_1, \ldots, z_5, v_1, \ldots, v_5]$, is
$I(V)_{D_4}=(v_5 v_2 v_4-z_1 z_3)$. So $D_4$ is directly
bipartite, $\{z_1, z_3\}$ is a source cover and $\{v_2, v_4,
v_5\}$ is a sink cover. The same result can be obtained by finding
directly sink and source covers of $D_4$ trough the vertex cover
of the undirected graph $K_{D_4}$.
\end{ex}

\end{document}